\documentclass{article}

\usepackage{arxiv}

\usepackage[utf8]{inputenc} 
\usepackage[T1]{fontenc}    
\usepackage{hyperref}       
\usepackage{url}            
\usepackage{booktabs}       
\usepackage{amsfonts}       
\usepackage{nicefrac}       
\usepackage{microtype}      
\usepackage{lipsum}		

\usepackage{graphicx}
\usepackage{natbib}
\usepackage{doi}
\usepackage{xcolor,colortbl}

\usepackage[linesnumbered,ruled,vlined]{algorithm2e}

\SetCommentSty{mycommfont}

\SetKwInput{KwInput}{Input}                
\SetKwInput{KwOutput}{Output} 
\SetKwInput{Kwset}{Set}

\usepackage{pgf-pie}
\usepackage{array,longtable}
\usepackage{xcolor,colortbl}
\usepackage{amsthm,amssymb,amsmath,mathtools}

\title{Initialization for Nonnegative Matrix Factorization: a Comprehensive Review}

\author{ \href{https://orcid.org/0000-0002-5731-8234}{\includegraphics[scale=0.06]{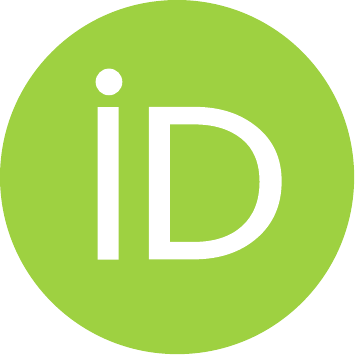}\hspace{1mm} Sajad Fathi Hafshejani}\\
	Department of Applied Mathematics\\
    Shiraz University of Technology\\
	Shiraz, Iran \\
	\texttt{s.fathi@sutech.ac.ir} \\
	\And
	{ Zahra Moaberfard} \\
	Department of Applied Mathematics\\
	Shiraz University of Technology\\
	Shiraz, Iran \\
	\texttt{zahra.moaberfard@gmail.com} \\
}

\renewcommand{\shorttitle}{\textit{arXiv} Template}

\hypersetup{
pdftitle={Initialization for Nonnegative Matrix Factorization: a Comprehensive Review},
pdfsubject={},
pdfauthor={S. Fathi Hafshejani,Z. Moaberfard},
pdfkeywords={Non-negative matrix factorization, Initialization algorithms},
}

\begin{document}
\maketitle

\begin{abstract}
	Non-negative matrix factorization (NMF) has become a popular method for representing meaningful
		data by extracting a non-negative basis feature from an observed non-negative data matrix.
		Some of the unique features of this method in identifying hidden data put this method amongst the powerful methods in the machine learning area. The NMF is a known non-convex optimization problem and the initial point has a significant effect on finding an efficient local solution. In this paper, we investigate the most popular initialization procedures proposed for NMF so far. We describe each method and present some of their advantages and disadvantages.  Finally,  some numerical results to illustrate the performance of each algorithm are presented.
\end{abstract}

	\textcolor[rgb]{0.31,0.24,0.57}{{\section{Introduction}}}
	Over the last few years, the  low-rank approximation,  which 
	is approximating a matrix by one whose rank is less than that of the original matrix has been an important technique and highly popular method in data science.  Low-rank approximations are fundamental and widely used tools for data analysis, dimensionality reduction, and data compression.  This method appears in many applications such as, image processing {\color{purple}\citep{Friedland}}, text data-latent semantic indexing, text mining {\color{purple}\citep{Elden,Skillicorn}}, and  machine-learning {\color{purple} \citep{Murphy,Kim}}. Low-rank approximations find two matrices of the much lower-rank that approximate a high-dimensional matrix $X$ such that:
	\begin{equation}
	X_{m\times n}\approx W_{m\times r}H_{r \times n},
	\end{equation}
	\noindent where
	\begin{equation}\label{rank}
	{r}<<\min(m,n)
	\end{equation}
	in which $r$ is so called {\it rank} of the matrix. There are very widespread matrix decompositions that give a low-rank approximation, for example we can refer to singular value decompositions (SVD) {\color{purple}\citep{Golub,8,17}}. It can provides the optimal rank and gives the appropriate of low-rank  approximation of a matrix {\color{purple}\citep{Datta,Sundarapandian,Trefethen}}. Unfortunately, these approximations usually do not actualize eligible structural constraints such as element-wise non-negativity {\color{purple}\citep{Fazel, Reinsel, Miller, Chu}}. For this reason, other concepts based on convex optimization have been developed such as URV, SDD, PCA, ICA, CUR, QR, and NMF {\color{purple} \citep{Petros,Smilde,Meyer,Tamara,23,14,18}}.
     Each of these approaches is based on different constraints that characterize the final properties of the matrix factors, leading to different optimization problems and numerical algorithms that must be used. 
	
	Non-negative Matrix Factorization (NMF) is an unsupervised data decomposition technique, akin to latent variable analysis, that can be used for feature learning, topics recovery, clustering, temporal segmentation, filtering, and source separation coding as with vector quantization. As a matter of fact, this method obtains parts-based, compression, and discriminant representation of the original data as well as enhancing the interpretability by using decomposes the main matrix into additive parts {\color{purple}{\citep{Daniel}}}. There have been some significant developments {\color{purple}{\citep{Aggarwal}}} in using  NMF for computation of the linear part-based representation of non-negative data. Therefore,  NMF can be considered as a method in the machine learning area which enhances the interpretability of the results. Therefore,  interpretability can be considered as one of the advantages of NMF. In addition, robustness is another property of NMF that can be applied to handle noise and estimates the missing values.
		Suppose that $X$ is a real and non-negative  $m \times n$ matrix.  NMF finds two real and non-negative matrices $W\in\Bbb{R}^{m \times r}$ and $H\in\Bbb{R}^{r\times n}$,
	such that:
\begin{equation}
	X\approx WH.\nonumber
\end{equation}
	  There are various approaches to find matrices $W$ and $H$.  An efficient way for this purpose is to apply optimization tools. To use optimization algorithms, we need a criterion for measuring the difference between the original matrix, that is, $X$, and the results, i.e, the matrix $WH$. This measure is called {\it objective function }  and can be considered as a measure for denoting an error. Here, we review two of the most common measures:
	\begin{itemize}
		\item Frobinous norm-based algorithm (SED): In this case, we use  the Euclidean distance between origin matrix $X$ and its approximation $WH$ as the similarity measure to derive the following objective function, which is based on the Frobenius norm for matrices:
		\begin{equation}\label{formmatrix1}
		\min_{W,H}f(W,H)= \frac{1}{2}\|X-WH\|^{2}_{F},~~\text{such that}~ W\geq 0 ~\text{and}~ H\geq 0.
		\end{equation}
		\item  Divergence-based algorithm (GKLD): This case is the most popular in real applications, the corresponding objective function that characterizes the similarity between matrix $X$ and matrix $WH$ (called divergence) and is given by {\color{purple}{{\citep{Lee01}}}}:
		\begin{equation}\label{formmatrix}
		\min_{W,H\geq 0}KL(X;WH)=\sum_{i,j}(X_{ij}\log \frac{X_{ij}}{(WH)_{ij}}-X_{ij}+(WH)_{ij}).
		\end{equation}
	\end{itemize}
Note that objective functions given by ({\color{purple}\ref{formmatrix1}}) and ({\color{purple}\ref{formmatrix}}) are non-convex in both $W$ and $H$. So, iterative approaches suggested for solving them guarantee to converge to some local minimum (more precisely, stationary points), but require initialization mechanisms that can greatly affect their convergence rate. ``Good" initial values for NMF are defined as follows {\citep{Boutsidis}}:
	\begin{itemize}
		\item One that leads to rapid error reduction and faster convergence.
		\item One that leads to better overall error at convergence.
	\end{itemize}
	There have been numerous results devoted to the initialization approaches for the  NMF {\citep{Casalino,math}}.  So, it seems that a systematic survey is of necessity and consequence. This review paper will summarize the most existing initialization strategies for NMF. We first present some common methods for solving NMF problem and then focus on the initialization approaches for this problem. We collect the most common initialization algorithms and investigate the advantages and disadvantages of them. Finally, we perform Lee's algorithm to compare the efficiency of the initialization methods. We perform the algorithm on the ORL  dataset consisting of face images and compare the results. 
	
	The rest of this paper is organized as follows. Section 2 briefly reviews the NMF problem and presents some common approaches for solving it and classifies the existing initialization methods.  In Section 3 we have a comprehensive review of random initialization seeding methods and present their algorithm.  We preset various clustering initialization strategies in Section 4.  Heuristic schemes for initialization are invested in Section5. In Section  6, we review some low-rank approximation methods.  In Section 7 we present some numerical results of performing Lee's Algorithm on the ORL database to demonstrate the performance of each initialization strategy. We finally
	end up the paper by giving some concluding remarks in Section 8.
	\textcolor[rgb]{0.31,0.24,0.57}{\section{NMF methods}}
	In this section, we review some common approaches for solving NMF. We start this section with recall multiplicative update rules, i.e.,  the SED-MU and GKLD-MU proposed by Lee and Seung  {\color{purple}\citep{Lee01}}. These methods have still been widely used as the baseline. The SED-MU update the matrices $W$ and $H$ by using the following strategy:
	\begin{eqnarray}
	W^{k+1}_{ia}&=&W^{k}_{ia}\frac{(X{H^{k}}^{T})_{ia}}{(W^{k}H^{k}{H^{k}}^{T})_{ia}},
	\qquad\forall i,a; \label{W_Update}\\
	H^{k+1}_{bj}&=&H^{k}_{bj}\frac{({W^{k+1}}^{T}X)_{bj}}{({W^{k+1}}^{T}W^{k+1}H^{k})_{bj}},
	\qquad\forall b,j.\label{H_Update}
	\end{eqnarray} 
	Moreover, the GKLD-MU can be formulated as:
	\begin{eqnarray}
	&&W_{ia}\leftarrow W_{ia}\sum_{j}\frac{X_{ij}}{(WH)_{ij}}H_{aj},\nonumber\\ &&W_{ia}\leftarrow \frac{W_{ia}}{\sum_{j}W_{ja}},\qquad ~~\text{and}
	~\quad H_{aj}\leftarrow H_{aj}\sum_{i}W_{ia}\frac{X_{ij}}{(WH)_{ij}}.\nonumber
	\end{eqnarray}
	It is proven that  the multiplicative update rules converge to a local
	minimum {\color{purple}\citep{Lee01}}.
	
	Another popular approach for solving NMF applying SED as the objective function  is called {\it Alternating Non-negative Least Squares} (ANLS), which is alternating least squares (ALS) modified under the non-negativity constraint.
	This approach finds two matrices $W$ and $H$ by solving the following optimization problems:
	\begin{eqnarray}
	&&W^{k+1}=\arg\min_{W\geq 0}f(W,H^K)=\frac{1}{2}\|X-WH^k\|_F^2,\label{f1}\\
	&&H^{k+1}=\arg\min_{H\geq 0}f(W^{k+1},H)=\frac{1}{2}\|X-W^{k+1}H\|_F^2.\label{f2}
	\end{eqnarray}
	Although the original problem ({\color{purple}\ref{formmatrix1}}) is non-convex and NP-hard with respect to variables $W$ and $H$, the sub-problems ({\color{purple}\ref{f1}}) and ({\color{purple}\ref{f2}}) are convex problems. However, these subproblems may have multiple optimal solutions because they are not strictly convex {\color{purple}\citep{g3}}.
	
	To accelerate the convergence rate, one popular method is to apply gradient descent algorithms with additive update rules. Other techniques such as conjugate gradient, projected gradient, interior point method, and more sophisticated second-order schemes like Newton and Quasi-Newton methods  are also in consideration {\color{purple}\citep{g1,g2,g3,g4,g5}}. To satisfy the non-negativity constraint, the updated matrices are brought back to the feasible region, namely the non-negative orthant, by additional projection, like simply setting all negative elements to zero.
	
 Now, we present a general framework for solving NMF in Algorithm {\color{purple}\ref{alg:Work_NMF}}, which  starts with the given non-negative matrix $X$ and the initial matrices, that is, $W^0$ and $H^0$. Then  it tries to find two non-negative matrices $W$ and $H$ such that the value of $\|X-WH\|$ is minimized. To do so, first algorithm checks the  stop condition. If the stop condition is not true, then two matrices $W$ and $H$ will be  updated with some common roles. The algorithm repeats this same process until the stop condition is true. In this case,  a appropriate solution, i.e., two matrix  $W$ and $H$ for origin problem is obtained.

\begin{algorithm}[H]
\DontPrintSemicolon
  
  \KwInput{$X,$ $W^{0},~H^{0},~\varepsilon$, and the  rank of approximation $r$}
  \KwOutput{$W$ and $H$}
  \Kwset{Objective function: $F=\|X-WH\|$}
  Checking Stop Condition(s)\;
  \eIf{$\|x^{r}-x^{r+1}\|<\varepsilon$ or finish iterations,}{
   finish algorithm
   }{
  	Trial step calculation: Update $W$ and $H$}
\caption{\textbf{The generic NMF Algorithm}}
\label{alg:Work_NMF}
	\end{algorithm}

	\begin{figure}[h!]
		\begin{center}
			\includegraphics[width=.75\textwidth]{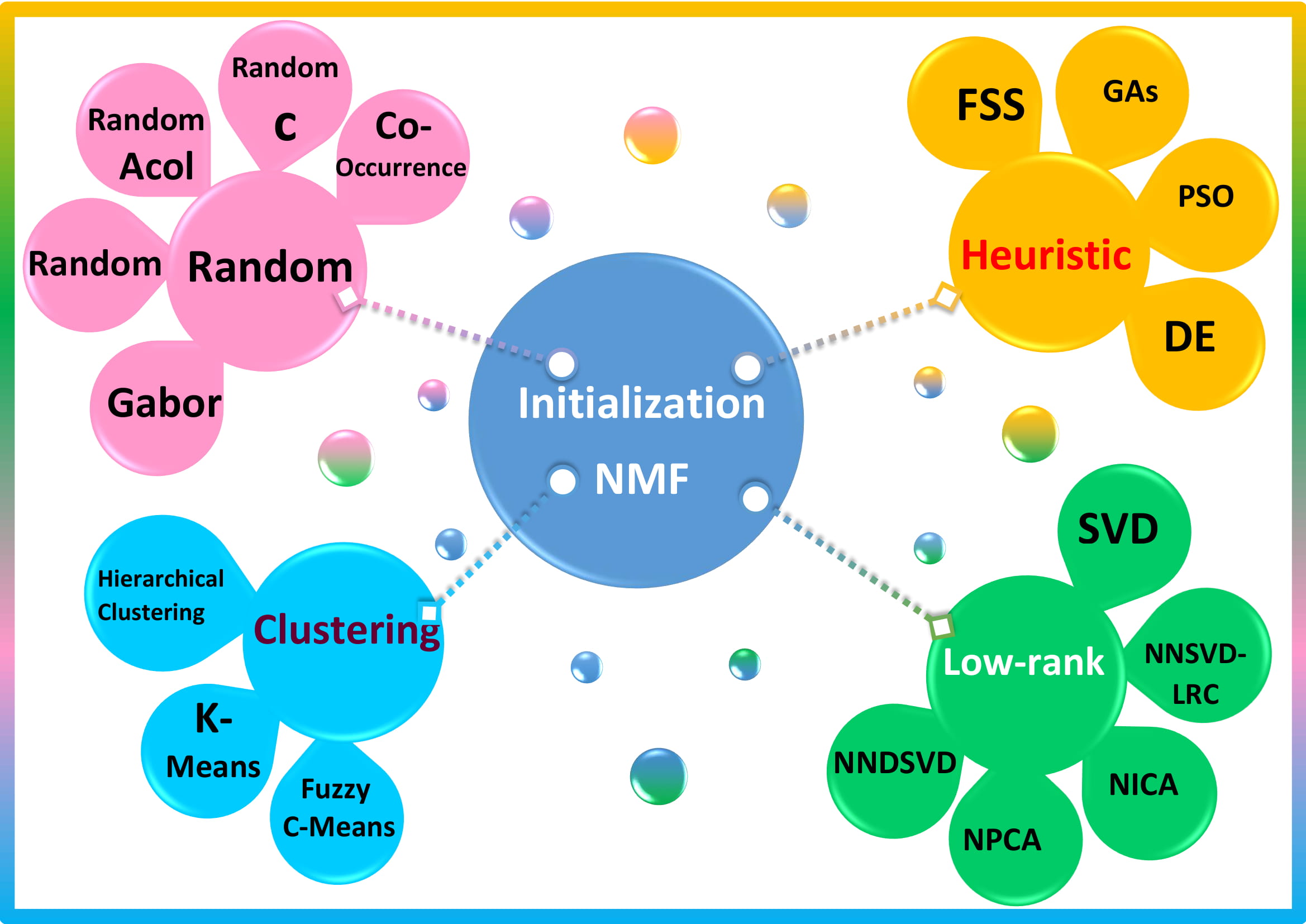}
			\caption{Proposed scheme for initialization NMF. }
			\label{fi}
		\end{center}
	\end{figure}
	As mentioned above, in iterative methods for solving NMF, the matrices $W$ and $H$ are obtained in such a way that the value of the objective function is minimized. Based on the non-convexity property for NMF, it generally does not guarantee a unique solution and its solution is dependent on choosing initialization for $W$ and $H$ demonstrated as
	$W^{0}$ and $H^{0}$ in this paper. A good choice for initializing can significantly affect the rate of convergence of the algorithm and considerably reduces the value of the cost function. Therefore, the goal of this paper is to investigate the initialization methods for NMF. The existing initialization approaches for NMF can be  classified  into four categories, as  shown in  Fig. {\color{purple}\ref{fi}}.  
	Random schemes, which only use the random strategy; Clustering schemes which profit the clustering strategy; Heuristic schemes, which are based on Population-Based Algorithms (PBAs); Low-rank Approximation-Based schemes, which works based on decreasing the matrix rank.
	
	Random strategy can be categorized into five sub-classes
	\begin{itemize}
	    \item Random, which suggests initial matrix by using random,
	    \item Random Acol, which calculates initial matrix $W$ by getting an average of  $q$ random columns of the matrix $X$
	    \item Random C, which calculates initial matrix $W$ by getting an average of the chooses q columns randomly from the longest (in the 2-norm) columns of $X$ {\color{purple}\citep{Casalino,Albright}}.
	    \item Co-Occurrence, which computes matrix $W$ by using $X^TX$ {\color{purple}\citep{Albright}}.
	    \item Gabor-based, which calculates the matrix $W$ by using Gabor wavelet and it is suitable for image datasets {\color{purple}\citep{Zheng07}}.
	\end{itemize}
	Correspondingly, Clustering strategy is categorized
	into three subclasses: 
	\begin{itemize}
	    \item  K-means, which use the K-means algorithm for initialization matrix $W$ {\color{purple}\citep{Xue08}}
	    \item Fuzzy C-means, which works based on the fuzzy roles {\color{purple}\citep{Zheng07,Alshabrawy,Rezaei}}
	    \item Hierarchical Clustering, which groups similar objects into groups called clusters {\color{purple}\citep{kim07}}. 
	\end{itemize}
	Besides, Heuristic Schemes  is categorized into four
	subclasses: 
	\begin{itemize}
	    \item Genetic Algorithm {\color{purple}\citep{Stadlthanner,Snasel,Price}}
	    \item Particle Swarm  Optimization
	    \item Differential Evolution
	    \item Fish School Search
	\end{itemize}
	Finally, Low-rank Approximation-Based is categorized into four subclasses
	\begin{itemize}
	    \item Singular Value Decomposition, which works based on SVD decomposition {\color{purple}\citep{Boutsidis}}
	    \item Nonnegative Singular Value Decomposition with Low-Rank Correction which generates a positive matrix {\color{purple}\citep{syed}}.
	    \item Non-negative PCA, which works based on PCA algorithm {\color{purple}\citep{Zheng07,Zhao08}}
	    \item Non-negative ICA, which works based on  ICA algorithm {\color{purple}\citep{Kitamura,Oja,Benachir}}
	\end{itemize}		
	In the following sections, we will discuss in detail each of these methods.


	\textcolor[rgb]{0.31,0.24,0.57}{\section{Random Schemes}}
	Random initialization is the benchmark used in the vast majority of NMF studies. Among several random initialization mechanisms, we selected five different random initialization strategies (which require low computational costs but have the drawback of generating poor informative initial matrices),
	namely Random, Random C, Random ACOL, Co-Occurrence, and Gabor-based initialization.  We describe them in the rest of this section. 
	
	\textcolor[rgb]{0.31,0.24,0.57}{\subsection{Random}}
	Probabilistic concepts can be used as an effective method for initializing the NMF. Over the past two decades, probabilistic approaches have been established to compute matrix approximations, forming the field of randomized numerical linear algebra {\color{purple}\citep{sandler}}. Random initialization is one of the common methods for the NMF algorithm that relies on a random selection of columns of the input matrix. However, the quality and reproducibility of the NMF result are rarely questioned when using random initialization.
	Different initialization of randomness and starting point will lead to different answers, so the algorithm should be run for several instances to select the best results of a local minimum. Random strategy can be used in many geometric initialization for NMF, in which the columns are selected in a more sophisticated way {\color{purple}\citep{pro1,pro2,pro3}}.
	
	In general, the randomized algorithms have shown their advantages for solving the linear least squares problem and low-rank matrix approximation. These methods have a low computational cost, but for some cases, the convergence rate to local minima and the qualitative solution is not guaranteed. Although randomness does not deliver reproducible results and does not generally provide a good first estimate for NMF algorithms.
	
	In the standard NMF algorithm, $W^0$ and $H^0$ are two non-negative matrices, where they have drawn from a uniform distribution, usually within the same range as the target matrices entries. This strategy is in-expensive and sometimes provides a good first estimation for the NMF algorithm. 	The first random initialization methods proposed by Lee and Soung {\color{purple}\citep{Lee01}}. Later on,
	this approach has been applied for various NMF algorithms, such as classical matrix factorization {\color{purple}\citep{Mahoney,Drineas,Casalino}}. In {\color{purple}\citep{Wangli}} Wang and Li proposed an algorithm working based on random projections to efficiently compute the NMF. Later, Tepper and Sapiro {\color{purple}\citep{Tepper}} in 2016 suggested a method that compressed the NMF algorithms based on the idea of bilateral random projections.
	While these compressed algorithms reduced the computational load considerably. However, this strategy is used in most of the NMF algorithms but has a drawback, that is, the algorithm needs multiple runs and in any performing, a different starting point is selected.  This significantly increases the computation time needed to obtain the desired factorization.
	To tackle this problem, several methods with different approaches for better seeking of NMF have been suggested, for example, computing a reasonable starting point from the target matrix itself.  Their goal is to produce deterministic algorithms that need to run only once, still giving meaningful results (e.g. Clustering, SVD) that in the following we will discuss.
	
	
	\textcolor[rgb]{0.31,0.24,0.57}{\subsection{Random Acol}}
	Random Acol forms an initialization of each column of the basis matrix $W^0$ by averaging $q$ random columns of matrix $X$ {\color{purple}\citep{Amy2}}. 
	Algorithm {\color{purple}\ref{alg:RandomAOLC}} presents a generic framework for initialization NMF based on  Random Acol.

	\begin{algorithm}[H]
\DontPrintSemicolon
  \KwInput{ $X, q, r$}
  \KwOutput{$W$}
  \Kwset{$k=0$}
 \While{$k\leq r$}{
 $k=k+1$\;
 Select $q$ columns  of matrix $X$ as random\;
 $s$=mean of $q$ columns\;
$W_{:,k}=s$
 }
\caption{\textbf{ Random ACOL Initialization algorithm}}
\label{alg:RandomAOLC}
	\end{algorithm}

Random Acol initialization builds basis vectors from the given data matrix; hence, as observed in {\color{purple}\citep{Albright}}, when the matrix $X$ is sparse, this initialization scheme forms a sparse initial basis matrix $W^{0}$, which represents a more reasonable choice compared to the random initialization. However, the performance of NMF algorithms initialized by Random Acol scheme is comparable with those of random initialization {\color{purple}\citep{Casalino}}. Nevertheless, Random Acol has one clear advantage over random initialization and it is creating a very sparse $W^{0}$, but this method is also very inexpensive but easy to implement.
	
	
	\textcolor[rgb]{0.31,0.24,0.57}{\subsection{Random C}}
	Random C initialization is similar to Random Acol initialization with only one main difference. In fact, it chooses \emph{q} columns randomly from the longest (in the 2-norm) columns	of the matrix $X$, which generally means the densest columns since our text matrices are so sparse. This method is also fairly inexpensive and easy to implement and is summarized in Algorithm {\color{purple}\ref{alg:RandomC}}.

	\begin{algorithm}[H]
\DontPrintSemicolon
  \KwInput{$X\in\Bbb{R}^{m\times n} $, $r$ and $q$}
  \KwOutput{$W$}
  \Kwset{$k=0$}
 \While{$k\leq r$}{
 $k=k+1$\;
 Find \emph{q} of the longest (in the 2-norm sense) columns of $X$\;
 $s$=mean of $q$ columns\;
$W_{:,k}=s$
 }
\caption{\textbf{Random C Initialization}}
\label{alg:RandomC}
\end{algorithm}

	It is shown that Random C initialization yields better results than the Random Acol initialization for either the asymmetric or the symmetric formulations in NMF {\color{purple}\citep{Casalino,Albright}}. Thus, the Random C initialization is more suitable compared to the Random Acol.  Despite having a low computational cost and providing a more realistic first estimate of the sources compared to random initialization, these methods suffer from a lack of reproducibility.
	
	
	\textcolor[rgb]{0.31,0.24,0.57}{\subsection{Co-Occurrence}}
Co-occurrence is a powerful tool for discovering the relationships	between heterogeneous collections of attributes or events.  Typically, if two such features frequently co-occur
throughout a database, it is assumed that they correspond to traits of the same object, concept, or process. The co-occurrence scheme first forms a term co-occurrence matrix $XX^{T}$. Next, this method randomly chooses the \emph{k} columns of the initial factor $W^{0}$ among the	densest columns of the co-occurrence matrix and generates $H^{0}$ (when required) via the random initialization  {\color{purple}\citep{sandler}}.
The co-occurrence scheme has the advantage of producing a
basis matrix that includes some hidden information on the initial data (i.e., term-term similarities when a document clustering scenario is considered). However, it requires a higher computational cost than simple random initialization. The co-occurrence method is very expensive for two reasons. First, if $m \gg n$, which means $C =XX^{T}$ is very large and often very dense too. Second,
the algorithm for finding $W^{0}$ is extremely expensive, making this method impractical. As evidenced by some
authors the Random C and co-occurrence initializations suffer from lack of diversity {\color{purple}\citep{Albright}}.
	
	
	\textcolor[rgb]{0.31,0.24,0.57}{\subsection{Gabor-based initialization}}
	Gabor wavelet is a powerful tools in image feature extraction  defined by {\color{purple}\citep{Zheng07}}:
	\begin{equation}
	\psi_{\mu,v}(z)=\frac{\|k_{\mu,v}\|^2}{\sigma^2}e^{(\frac{-\|k_{\mu,v}\|^2\|z\|^2}{2\sigma^2})}\left[e^{-ik_{\mu,v}z}-e^{-\frac{\sigma^2}{2}}\right]
	\end{equation}
	in which  $\mu$ and $v$ denote the orientation and scale of the Gabor kernels and the wave vector $k_{\mu,v}$ is given by:
	\begin{equation}
	k_{\mu,v}=k_ve^{-i\phi_{\mu} }
	\end{equation}
	So, the Gabor feature representation of an image $I(z)$ is obtained by:
	\begin{equation}
	G_{\mu,v}(z)=I(z)*\psi_{\mu,v}(z)
	\end{equation}
	where $z=(x,y)$ and $*$ is the convolution operator.
	When an image convolves with Gabor wavelets, the image is transformed into a set of image features at certain scales and orientations. Therefore, the image can be reconstructed from these image features. Motivated by this point,  Zheng et al. {\color{purple}\citep{Zheng07}} applied  the Gabor-based method to initialize NMF. The advantage of this method is that it is very suitable for image data sets.


	\textcolor[rgb]{0.31,0.24,0.57}{\section{Clustering Schemes}}
The clustering-based method is one of the common approaches in initialization for  NMF. Since this method produces a	summarized view of data helping the analyst to understand data by means of compact and informative representations of large collections	of samples {\color{purple}\citep{Berthold}}.	{The NMF as a clustering method can be traced back to work by Lee and Seung {\color{purple}\citep{Lee01}}. But, the first work that explicitly demonstrates that it was done by Xu et al. {\color{purple}\citep{Xu}}}. Typically, clustering algorithms are initialized by random strategy. Moreover, these methods have good results in environmental research in public health {\color{purple}\citep{Chretien}}, signal and image	processing  {\color{purple}\citep{Cichocki4}}. If  NMF method is	considered as a clustering process, the initialization strategy can	be obtained based on the results of clustering algorithms and fuzzy clustering. {There are various  types of
	clustering approaches, for example, supervised/unsupervised, hierarchical/partitional,  hard/soft, and one-way/many-way (two-way clustering is known as co-clustering or bi-clustering) among others.} Clustering-based initialization schemes will provide more realistic source estimates compared to low-rank approximation methods, but they can be computationally expensive. Furthermore, clustering	methods usually require some initialization themselves. Most of the proposed	initialization methods have been compared with random initialization in terms	of convergence rate and/or quality of the solution. However, different random initializations will lead to different NMF results, making it a questionable	reference. It is unclear how previous studies have dealt with the lack of reproducibility.	In this case that prototype-based clustering is a convenient method for the problem at hand, NMF could be	a valid tool. NMF has been widely used in clustering applications {\color{purple}\citep{Shahnaz, Xu}},	where the factors $W$ and $H$ have been interpreted in terms	of cluster centroid and cluster membership, respectively. On the other hand, the divergence-based NMF algorithm is not utilization {\color{purple}\citep{Wild04,Wild02}}. There are several initialization methods that work based on a clustering scheme. Most of these methods have used the Euclidean distance between the input matrix and the NMF approximation.

Many different clustering methods exist in the literature, such as	hierarchical clustering, prototype-based clustering,
and density-based clustering. Hierarchical clustering yields a	collection of nests groups of data, while in prototype-based clustering groups are represented in a compressed form through a prototype, i.e., an element belonging to the same domain of data. In density-based clustering, groups	are formed in regions of data space where data are more crowded. The choice of the	most appropriate method is up to the data analyst.	In the following, we will concentrate on three well-known	clustering schemes i.e., K-means, Fuzzy C-means,  and  Hierarchical Clustering.

	
	\textcolor[rgb]{0.31,0.24,0.57}{\subsection{K-means}}
	The K-mean method is a clustering technique used to grouped
	similar patterns in given features.
	{The K-means (Grst introduced it in 1960 {\color{purple}\citep{Forgy}})
		is the most widely used clustering technique {\color{purple}\citep{Hartigan}}}.
	This method represents points in the $k$-space that are the centers of
	clusters of nodes with the characteristic that they minimize the sum of squared distance deviations of
	the points in each cluster from
	the assigned cluster "centroid". The K-means algorithm is an
	iterative algorithm for minimizing the sum of distance between
	each data point and its cluster center (centroid) and tries
	to minimize the sum-squared-error criterion.
	Generally, K-means method seeks to partition the data set $X$
	into $k$ disjoint clusters so that each
	point in the cluster is ``closer'' to the centroid associated
	with that cluster than it is to the other $k-1$ centroids in
	the Euclidean sense.
	
	As the K-means factor is added to NMF, it gives prominent
	importance in clustering with extracted features.

	The theoretical connection between factorization NMF with additional orthogonal constraints on its factors, K-means, and spectral clustering was demonstrated in {\color{purple}\citep{Ding}}. While the mathematical
	equivalence between orthogonal NMF and a weighted variant of spherical
	K-means was proved together with some indications about the cases in
	which orthogonal NMF should be preferred over K-means and spherical K-means.
	
	The objective functions for K-means is defined as:
	\begin{equation}\label{k-means}
	M=\sum^{l}_{j=1}\sum^{n}_{i=1}\|x_{i}^{(j)}-c_{j}\|^{2},
	\end{equation}
	where $x$ is the feature vector, $c_j$ denotes the center of the
	cluster and $j$ is the number of the cluster centers.
	The theoretical connection between K-means and NMF can be presented as:
	\begin{equation}\label{k-nmf}
	\min\sum^{n}_{i=1}\|x_{i}-w_{\sigma_{i}}\|^{2}_{2}=\min\|X-WH\|^{2}_{F},
	\end{equation}
	where  $W$ and $H$ are two non-negative  matrices. Moreover, $x_{1},...,x_{n}$ denote the columns of $X$ and $w_{1},...,w_{k}$ are the ${k}$ centroids and $\sigma_{i}=j$ when \emph{i}-th point is assigned to ${j}$-th cluster ($j\in{1,...,k}$). {Algorithms find minimum of  $M$  often apply iterative gradient descent approaches. These algorithms usually converge to local minima.}
	
	There are several methods for initializing the NMF based K-mean. We point out some of them here.
	\begin{itemize}
		\item The initial basis matrix W is constructed by using the K-means clustering approach and the initial matrix H is considered as a random matrix.
		\item The initial basis matrix $W$ is constructed by using the K-means clustering strategy and the initial matrix $H$ is calculated by $H=W^TX$ and then the absolute value function is used for all elements in $H$ in order to satisfy the initial constraint of NMF. 
		\item The initial basis matrix $W$ is obtained by 
		using the cluster centroids obtained from K-means 
		clustering. The initial matrix $H$ is obtained by 
		$H=W^TX$ and then all negative 
		elements in $H$ are transferred to zero in order to satisfy the initial 
		the constraint of NMF.
		\item The initial basis matrix $W$ is obtained by using the cluster centroids obtained from K-means clustering.  The value of the membership degrees of each data point is calculated by:
		\begin{equation*}
		h_{kq}=\frac{1}{\left(\sum_{k'}^{k}(\frac{d(x_q,c_{k'})}{d(x_q,c_k)})\right)^{\frac{2}{1-m}}}
		\end{equation*}
		where $d(.)$ denotes the Euclidean distance between the two points, $x_q$  represents the $q$-th data point and $c_k$ represents the $k$-th cluster centroid. Moreover,  the fuzzification parameter is denoted by $m$. The initial matrix 
		$H$ is then obtained by using the membership degrees above.	
		
	\end{itemize}

	Many random initialization methods for the K-means algorithm have been proposed so far. Most classical methods are random seed {\color{purple}\citep{Forgy,Anderberg}} and random partition {\color{purple} \citep{Anderberg}}. Random seeds randomly select \emph{k} instances (seed points) and assign each of the other instances to the cluster with the nearest seed point. Random partition assigns each data instance into one of the \emph{k} clusters randomly. To escape from getting stuck at a local minimum, one can apply \emph{r} random starts.
	
	To improve the performance of divergence-based NMF algorithm that works based on Xue's idea {\color{purple}\citep{Xue08}}, a new method using the K-means and combination of normalizing technique with set divergence as the similarity measure in clustering to find the base vectors for NMF initialization and search of the Centroids was first proposed by {\color{purple}\citep{Xue08}}. The authors used the $L_1$ norm to normalize their algorithm. The proposed algorithm primary ones to utilize NMF and it stops when the number of clusters does not change.

	\textcolor[rgb]{0.31,0.24,0.57}{\subsection{Fuzzy C-means}}
	The fuzzy set theory introduced by Zadeh 1965 provides a
	powerful analytical tool for the soft clustering method. The
	Fuzzy C-Means is the best-known
	approach for fuzzy clustering, based on optimizing an objective function.  This concept has many applications as a convenient tool in clustering and has the most perfect algorithm theory. The FCM clustering algorithm can be considered as a variation and an extension for the traditional K-means clustering algorithm, in which for each data point a degree of membership or membership function of clusters is assigned.  It is proven the fuzzy clustering is an adaptation to noisy data and classes that are not well separated.  By considering this property of fuzzy clustering,  some research papers were done in this area. For example, Zheng and et al.  {\color{purple}\citep{Zheng07}} proposed the FCM concept. They used their strategy to initialization of NMF.  Another work in this area was done by Alshabrawy et al.  {\color{purple}\citep{Alshabrawy}}. They applied the FCM clustering technique to estimate the mixing matrix and to reduce the requirement for the sparsity of the Semi-NMF.  In another work, Rezaei et al. {\color{purple}\citep{Rezaei}} applied FCM to initialize NMF as an efficient method to  enhance NMF performance.

	\textcolor[rgb]{0.31,0.24,0.57}{\subsection{Hierarchical Clustering}}
	This method is motivated by common sense on ``part",
which is the smallest unit that has some perceptual meaning. 	For example, a face image consists of various parts, including eyes, nose, 	eyebrows, cheek, lip, and so on. Metaphorically, a pixel corresponds to an atom, and then, a part can be considered as a molecule. As atoms in a molecule perform a chemical reaction together, pixels that build a part should be grouped together. They introduced a ``closeness to rank-one" (CRO) measure in order to investigate whether row vectors in the sub-matrix show similar patterns or not.
The CRO measure is defined by:
	\begin{equation}\label{cro1}
	CRO(X_{(i,j,...,k),:})=\frac{\sigma_1^2}{\sum_{i=1}^{r}\sigma_i^2}=\frac{\sigma_1^2}{\|X_{(i,j,...,k),:}\|_F^2},
	\end{equation}
	where $\sigma_1\geq \sigma_2\geq ...\geq \sigma_r\geq 0$ are singular values of the sub-matrix $X_{(i,j,...,k),:}$ and $r$ denotes the rank of the $X_{(i,j,...,k),:}$.
	Algorithm {\color{purple}\ref{CRO}} presents a general framework for the CRO algorithms. It was first proposed in {\color{purple}\citep{kim07}}. They used this method to initialize the NMF algorithm. They compared their results with random initialization strategy and investigated how goodness-of-fit (GOF) and sparseness changes after the convergence of the standard NMF algorithm starting from these two different initialization methods.

	\begin{algorithm}[H]
\DontPrintSemicolon
  \KwInput{$X\in\Bbb{R}^{m\times n} $ with $X\geq 0 $ and $r<<\min{(m,n)}$}
  {Assign each row vector $X_{1,:},X_{2,:},...,X_{m,:}$ into each own clusters $C_1,C_2,...,C_m$}\;
 {Calculate CRO in (\ref{cro1}) between every pairs of clusters}\;
 \While{$n$ clusters remains}{
 Find a pair of clusters with the largest CRO\;
 Merge them into a single cluster\;
 Compute CRO between the newly-merged cluster
				and remaining clusters
 }
		\caption{\textbf{CRO-based hierarchical clustering}}
		\label{CRO}
	\end{algorithm}

	\textcolor[rgb]{0.31,0.24,0.57}{\section{Heuristic Schemes}}
	
	An important aspect which has not been deeply investigated yet
	 is a proper initialization of the NMF factors in order
	to achieve a faster error reduction {\color{purple}\citep{Dong}}. Thus, several heuristics algorithms have been
	proposed to solve NMF problem. However, only a few studies combined NMF
	and Population Based Algorithms (PBAs) and both of them are based on population-based optimization algorithms.
	{\color{purple}\citep{Goldberg}} presented Genetic Algorithms (GA) which are global search
	heuristics that operate on a population of solutions using techniques encouraged from
	evolutionary processes such as mutation, crossover, and selection.
	{\color{purple}\citep{Stadlthanner}} investigated the application of GAs on sparse NMF
	for microarray analysis, while {\color{purple}\citep{asel}} proposed GAs for Boolean
	matrix factorization, a variant of NMF for binary data based on Boolean algebra. The
	results in these two papers are promising but barely connected to the initialization
	techniques introduced in this paper.
	In particle Swarm Optimization (PSO) {\color{purple}\citep{Kennedy}} each particle in
	the swarm adjusts its position in the search space based on the best position it has
	found so far as well as the position of the known best-fit particle of the entire swarm.
	In Differential Evolution (DE) {\color{purple}\citep{Price}} a particle is moved around
	in the search-space using simple mathematical formulation if the new position is an
	improvement the position of the particle is updated, otherwise, the new position is discarded.
	
	Algorithm {\color{purple}\ref{pbc}} presents a pseudo code for NMF initialization using PBAs which was first proposed in {\color{purple}\citep{Janecek}}. They used their strategy to initialize the NMF and compared obtained results with some other algorithms such as random, NNDSVD and showed that their method has better results in term of convergence.

\begin{algorithm}[H]
\DontPrintSemicolon
  \KwInput{$X\in\Bbb{R}^{m\times n} $ with $X\geq 0 $ and $r<<\min{(m,n)}$}
  \KwOutput{$W$ and $H$}
  \Kwset{$H^0=rand(r,n)$, $i=0$}
 \While{$i\leq m$}{
 $i=i+1$\;
 Use PBAs to find $w_i^r$ that minimizes $\|x_i^r-w_i^rH^0\|_F$\;
 set $W(i,:)=w_i^r$
 }
  \Kwset{$j=0$}
  \While{$j\leq n$}{
 $j=j+1$\;
 Use PBAs to find $h_j^c$ that minimizes $\|x_j^c-Wh_j^c\|_F$\;
 set $ H(:,j)=h_j^c $
 }

	\caption{\textbf{Pseudo Code for NMF Initialization using PBAs}}
		\label{pbc}
	\end{algorithm}

	There are two other papers that combine NMF and PBAs. In fact,  both of them are based on GAs. {\color{purple}\citep{Stadlthanner}} have investigated the application of GAs on sparse NMF for microarray analysis, while {\color{purple}\citep{Snasel}} have applied GAs for boolean matrix factorization,
	a variant of NMF for binary data based on Boolean algebra.

	\textcolor[rgb]{0.31,0.24,0.57}{\section{Low-rank   Approximation-Based}}
	In this section, we focus on the most important the Low-Rank (LR) approximation  algorithms. Initialization schemes based on LR decomposition strategies do not require a randomization stet.  The LR methods include strategies using the Singular
	Value Decomposition (SVD), Nonnegative Singular Value Decomposition with Low-Rank Correction (NNDSVD-LRC),   Non-negative Principal Component Analysis (NPCA), and   Non-negative Independent Component Analysis (NICA).

	
	\textcolor[rgb]{0.31,0.24,0.57}{\subsection{ Singular
			Value Decomposition}}
	
	The potential impact of the NMF and its extensions on scientific advancements might be as great as the other popular matrix factorization technique such as SVD, that is based on low-rank approximations. LR approximation using SVD has many applications over a wide spectrum of disciplines.
	For example,  image compression, similarly, text data latent semantic indexing {\color{purple}\citep{S.Deerwester}}, event detection in streaming data, visualization of a document corpus and etc.
	In particular, Boutsidis and Gallopoulos {\color{purple}\citep{Boutsidis}}, pointed out
	SVD-NMF has good properties under these two conditions:
	\begin{description}
		\item[$\bullet$] One that leads to rapid error reduction and faster convergence.
		\item[$\bullet$] One that leads to the overall error at convergence.
	\end{description}
	There exists a factorization with the following form:
	\begin{equation}\label{SVD1}
	X=W\Sigma H^{T}
	\end{equation}
	Let us denote orthogonal matrices as $W=[w_{1},w_{2},...]$ and $H=[h_{1},h_{2},...]$
	that include  the  left  and  right  singular  vectors  of $X$, respectively. Moreover,
	the matrix $\Sigma=diag(\sigma_{1},\sigma_{2},...,\sigma_{r})$  is a diagonal matrix containing the first ${r}$ singular values of $X$ and $(\sigma_{1}\geq\sigma_{2}\geq....\geq\sigma_{r}>0)$. The truncated SVD
	is the best rank-${r}$ approximate of matrix $X$, in either spectral norm or Frobenius
	norm {\color{purple}\citep{Eckart}}.
	In particular, the singular values decrease quickly as $\emph{i}$ increasing in most
	instances {\color{purple}\citep{Lijie}}, which means that some of the first singular values
	can contain almost all singular information of input matrix.\\

\begin{algorithm}[H]
\DontPrintSemicolon
  
  \KwInput{ $X\in\Bbb{R}^{m\times n} $ with $X\geq 0, $ and the rank matrix $r$.}
  \KwOutput{$W$ and $H$}
\textbf{Calculate $[u~ s ~v ~r] = truncated-SVD(X,r)$}\;
		\textbf{Calculate $W = abs(u(:,1:r))$}\;
		\textbf{Calculate $H = abs(s(1:r,:)*v)$}\;
		\caption{\textbf{Singular Value Decomposition Initialization}}
		\label{SVDal}
	\end{algorithm}

	This property allows us to use it to compress the matrix data by eliminating the small
	singular values or the higher ranks.
	From ({\color{purple}\ref{SVD1}}), the sum of all ${r}$ non-zero diagonal
	entries of the number of singular values are selected as:
	\begin{equation}\label{num_sing}
	\frac{sum_{i}}{sum_{k}}\geq 90 \;\;\; and \;\;\; \frac{sum_{i}-1}{sum_{k}} <90
	\end{equation}
	In {\color{purple}\citep{Boutsidis}}, the authors used an SVD-based initialization and showed anecdotical examples of speed up in the reduction of the cost function.  We present a generic algorithm for initializing NMF based SVD in Algorithm {\color{purple}\ref{SVDal}}. As mentioned before, SVD does not necessarily produce the non-negative matrices. So, some algorithms in this area change negative elements to positive or zero.

	This method have some drawbacks as:
	\begin{itemize}
		\item The interpretability of the transformed features. The resulting orthogonal matrix factors generated by the approximation usually do not allow for direct interpretations in terms of the original features because they contain positive and negative coefficients {\color{purple}\citep{Zheng07}}. 
		\item  This method suffers from the fact that the approximation error $\|X -W H\|_F^2$ of the initial factors $(W, H)$ increases as the rank increases which is not a desirable property for NMF initializations. 
	\end{itemize}
	Non-negative Double Singular Value Decomposition (NNDSVD) {\color{purple}\citep{Boutsidis}} is another method designed to enhance the initialization stage of the NMF.  This method contains no
	randomization and is based on two SVD processes, one approximating the data matrix, the other approximating positive sections of the resulting partial SVD factors utilizing an algebraic property of unit rank matrices. 
	NNDSVD can readily be combined with existing NMF algorithms. This property leads to the NNDSVD being considered as an efficient method for initializing the NMF.

	
	\textcolor[rgb]{0.31,0.24,0.57}{\subsection{Nonnegative Singular Value Decomposition with Low-Rank Correction }}
	The Nonnegative Singular Value Decomposition with Low-Rank Correction (NNSVD-LRC) was first proposed in {\color{purple}\citep{syed}}. This method works based on the SVD but it has some useful properties such as:
	\begin{itemize}
		\item  This method	generates sparse factors which not only provide storage efficiency but also provide better partbased representations and resilience to noise.
		\item It only requires a truncated SVD
		of rank $[\frac{r}{2}+1]$. 
	\end{itemize} 
	Here, we describe the NNSVD-LCR framework in Algorithm {\color{purple}\ref{NNSVD-LCR}} which was first proposed in {\color{purple} \citep{syed}}:

\begin{algorithm}[H]
\DontPrintSemicolon
  
  \KwInput{$X\in\Bbb{R}^{m\times n} $ with $X\geq 0 $ and $r<<\min{(m,n)}$}
  \KwOutput{$W$ and $H$}
  \Kwset{$p=\lfloor\frac{r}{2}+1\rfloor$}
 $[U, \Sigma, V ]$ = truncated-SVD$(X,p)$\\
 $Y_p=U\Sigma^{\frac{1}{2}}$, $Z_p=\Sigma^{\frac{1}{2}}V^T$\\
 $W(:,1)=|Y_p(:,1)|$, $H(1,:)=|Z_p(1,:)|$\\
 $i=2$; $j=2$\\
 \While{$i\leq r$}{
  \eIf{$i$ is even then}{
   $ W(:,i)=\max(Y_p(:,j),0)=\max(Z_p(:,j),0)$
   }{
   j=j+1\;
  $W(:,i)=\max(-Y_p(:,j),0)=\max(-Z_p(:,j),0)$}
$i=i+1$ }
 \caption{Nonnegative Singular Value Decomposition with Low-Rank Correction (NNSVD-LRC)}
 	\label{NNSVD-LCR}
\end{algorithm}

	\textcolor[rgb]{0.31,0.24,0.57}{\subsection{ Non-negative PCA}}
	Principal Component Analysis (PCA) is one of the best-known unsupervised feature extraction methods because of its conceptual simplicity and the existence of efficient algorithms that can implement it. Particularly in the face representation task, faces can be economically represented along with the eigenface coordinate space, and approximately reconstructed using just a small collection of eigenfaces and their corresponding projections (coefficients). It is an optimal representation in the sense of mean-square error. However, it presents some drawbacks (such as the presence of mixed sign values), and several research papers demonstrated that it outperforms NMF in many applications such as face recognition {\color{purple}\citep{Cichocki3,Guillamet}}.\\

	Principal components sequentially capture the maximum variability of $Y$ thus guaranteeing minimal information loss, and they are mutually uncorrelated. To clarify, consider the problem of human face recognition, where PCA has been largely adopted to obtain a set of basic images the eigenfaces, that can be linearly combined to reconstruct images in the original dataset of face {\color{purple}\citep{Turk}}.
	Here we describe this method. 
	Given the $m\times n$  matrix $X$ as that of in NMF, we define the  average vector $\psi$ as:
	\begin{equation*}
	\psi=\frac{1}{n}\sum_{i=1}^{n}X_i
	\end{equation*}
	The centered matrix can be calculated as:
	\begin{equation*}
	\bar{X}=(X_1-\psi,...,X_n-\psi).
	\end{equation*}
	Using SVD to compute the eigenvectors of the $\bar{X}^T\bar{X}$. The eigenvector
	matrix $W$ is constructed by keeping only $r$ eigenvectors
	(corresponding to the $r$ largest eigenvalues $\lambda_i$) as column
	vectors and $H$ is an $r\times n$ matrix containing the encoding
	coefficients. Moreover, the criterion of selecting $r$ is usually as follows
	\begin{equation*}
	\frac{\sum_{i=1}^{r}\lambda_i}{\sum_{i=1}^{n}\lambda_i}\geq \alpha,
	\end{equation*}
	in which $\alpha =0.9$. As we know that the NMF seeks to finds two non-negative matrices for initialization. So, the following non-linear operator can be used to transfer negative elements to zero.
	\begin{eqnarray}
	p(W)=p(W+)=\max(W_{ij},0)~\quad p(H)=p(H+)=\max(H_{ij},0).
	\end{eqnarray}
	
	There are several works in this area, for example, Zheng et al. {\color{purple}\citep{Zheng07}} proposed PCA-based 
	initialization method and, after obtaining $W$ and $H$, all 
	negative elements in these two matrices change to zero. In another work, Zhao et al. {\color{purple}\citep{Zhao08}}
	used the absolute value for all elements in the initial 
	matrices $W$ and $H$ after PCA initialization. With these initialization 
	methods, enhanced convergence rates as well as better 
	accuracy were achieved. 
	Geng et al. {\color{purple}\citep{Xiu-rui-16}} pointed out the NMF is  sensitive to noise (outliers), and used PCA to initialize the NMF.

	Despite the popularity of PCA, this method has two key drawbacks:
	\begin{itemize}
		\item It lacks sparseness
		(i.e., factor loadings are linear combinations of all the input variables), yet sparse
		representations are generally desirable since they aid human understanding (e.g., with
		gene expression data), reduce computational costs and promote better generalization in
		learning algorithms
		{\color{purple}\citep{Shen,Zass}}
		\item PCA is computationally expensive, the size of the
		covariance matrix is proportional to the dimension of the data. As a result, the computation of
		the eigenvectors and eigenvalues might be impractical for high-dimensional data.
	\end{itemize}
	
	\textcolor[rgb]{0.31,0.24,0.57}{\subsection{Non-negative ICA}}
	
	Independent component analysis (ICA) is another mechanism that used to extract a set of statistically independent source variables from a collection of mixed signals without having information about the data source signals or the combination process. The initialization methods using PCA or SVD are based on the orthogonality between the bases representing the data matrix $X$. However, it has been shown that the optimal NMF bases are along the edges of a convex polyhedral cone, which is defined by the observed points in $X$, in an M-dimensional space {\color{purple}\citep{29,30}}. Therefore, PCA and SVD may not be the best methods for the initialization in NMF. To avoid of this situation, some researchers proposed the utilization of NICA bases and estimated independent sources as the initial values of the basis and weight matrices, respectively {\color{purple}\citep{Kitamura,Oja}}. The numerical results provided faster and deeper convergence of the NMF cost function than the conventional methods.
	Benachir et al in {\color{purple}\citep{Benachir}} modified standard ICA taking
	into account the sum-to-one constraint and then eliminated
	some indeterminacies related to ICA using different strategies.
	Then, they used the outputs to initialize an NMF
	method.

	\textcolor[rgb]{0.31,0.24,0.57}{\section{Numerical results}}
	In this section, we present some numerical experiments of performing Lee's
	algorithm on the ORL dataset with some different initialization strategies.  The goal of this section is to compare the  accuracy of the algorithm based on different initialization approaches. In all experiments, we have used the stopping condition as:
	\begin{equation*}
	\|W^kH^k- W^{k-1}H^{k-1}\|\leq\epsilon,
	\end{equation*}
	where $\epsilon=10^{-10}$.
	
\begin{figure}[h!]
		\begin{center}
			\includegraphics[width=1\textwidth]{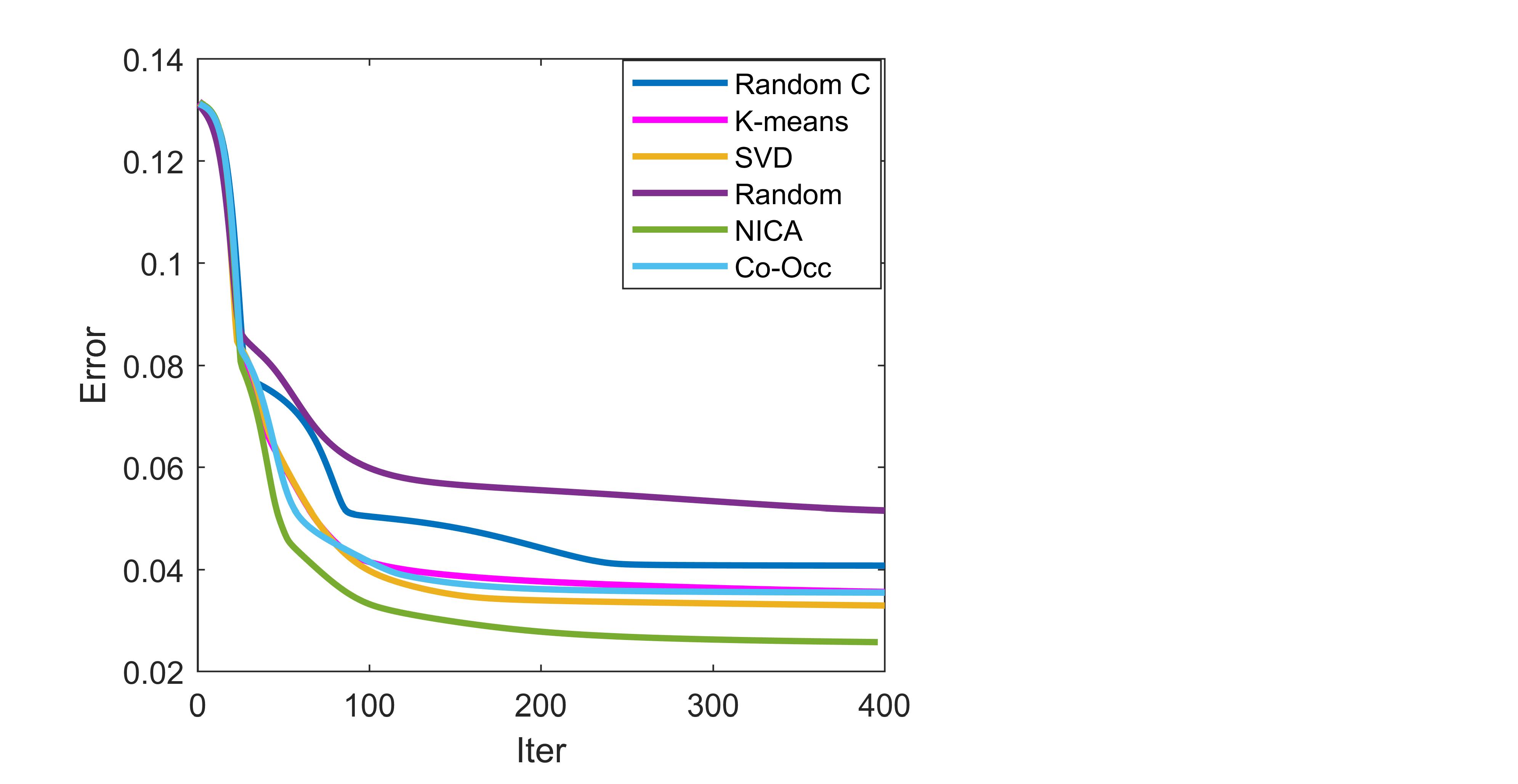}
			\caption{Number of training data is $320$.}
			\label{n8}
		\end{center}
	\end{figure}

\begin{figure}[h!]	
		\begin{center}
			\includegraphics[width=1\textwidth]{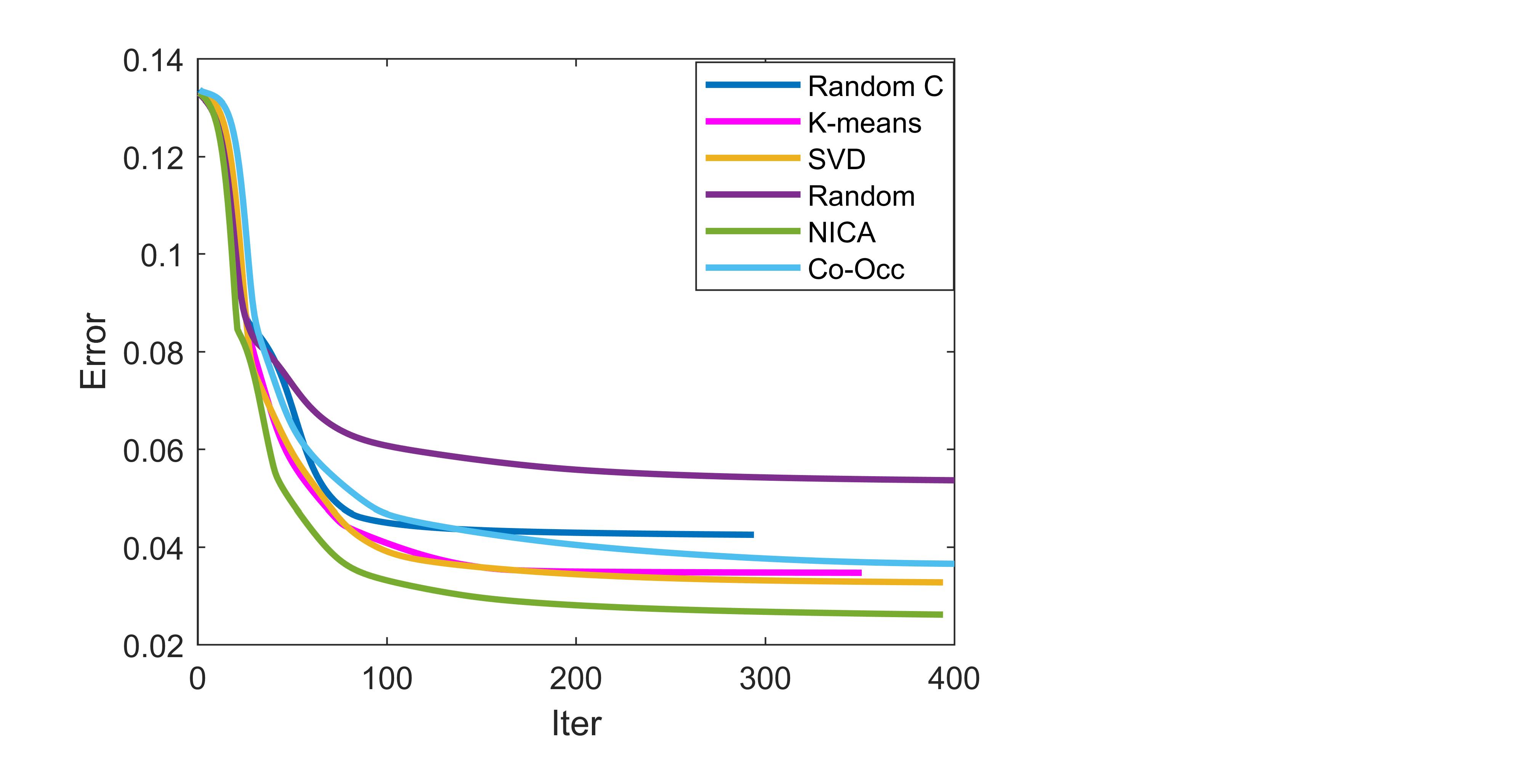}
			\caption{Number of training data is $200$.}
			\label{n5}
		\end{center}
	\end{figure}

\begin{figure}[h!]
		\begin{center}
			\includegraphics[width=1\textwidth]{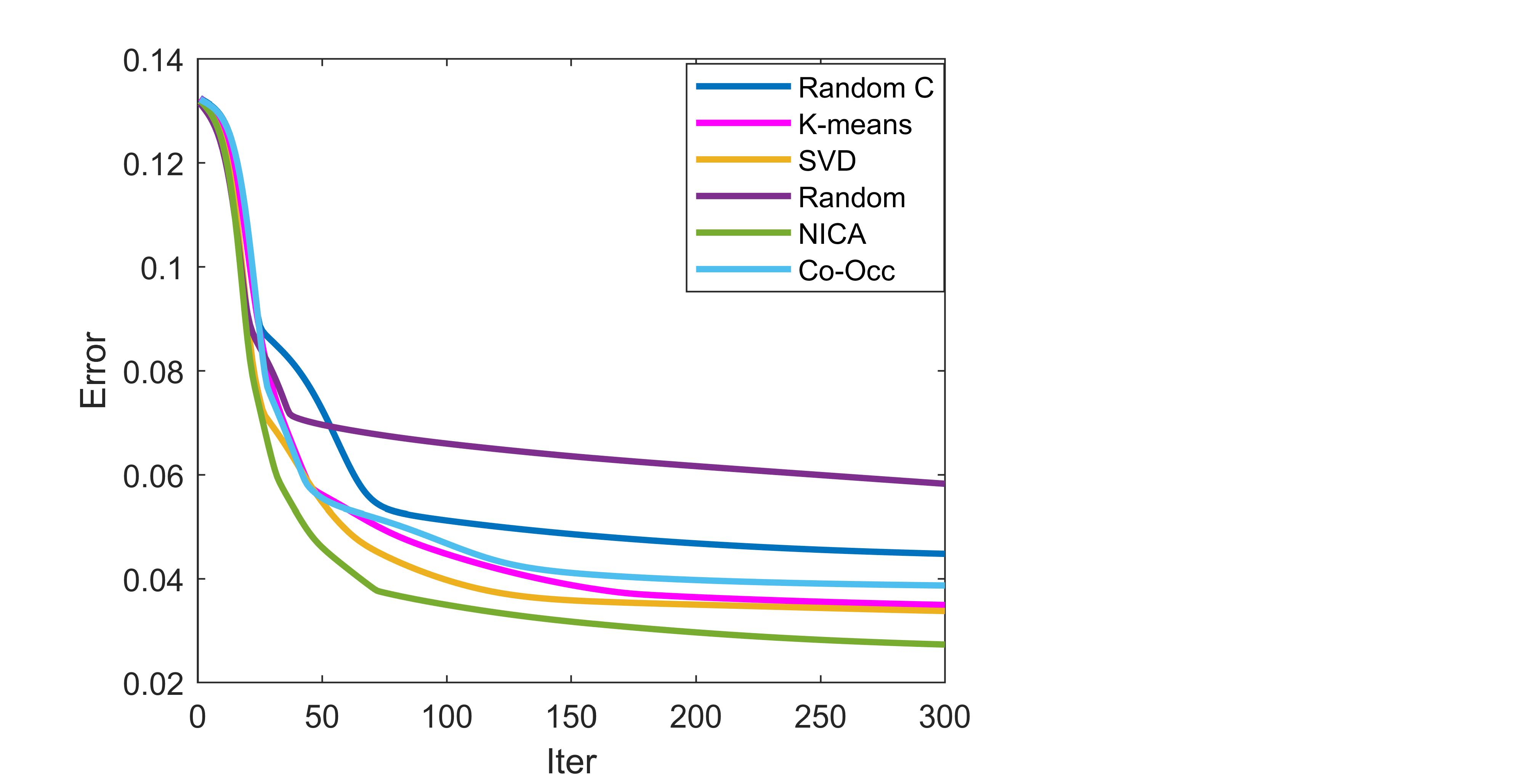}
			\caption{Number of training data is $120$.}
			\label{n3}
		\end{center}
	\end{figure}

\textcolor[rgb]{0.31,0.24,0.57}{{\subsection{Datasets and settings}}}
	All codes of the computer procedures are written in  MATLAB 2017 environment and are carried out
	on a PC (CPU 2.60 GHz, 16G memory) with the Windows 10 operation system environment. We initialize the NMF by using Random, Co-Occurrence, Random C, SVD, K-means, and NPCA strategies.   
	We perform the algorithm on the dataset ORL which has 400 images of 40 different classes that each of them has 10 images.  As we know that the number of images for train and test is very important in machine learning area. So, we run the algorithm with three different cases in terms of the number of data for training and testing and we report the accuracy of the algorithm in each case. In fact, we perform the algorithm with number of training as $\{320,200,120\}$.
	The error for all experiments is  calculated by:
	\begin{equation}
	Error:=\frac{\|X-WH\|_F}{\|X\|_F}.
	\end{equation}
	For random strategies, we performed  the algorithm 10 times and demonstrated the average results for each case.
	Fig. \ref{n8}~ shows approximation error for the case where the number of training data is $320$. As we see that, the NICA has the best result than the other algorithms. Moreover, the random strategy has the highest error. Fig. \ref{n5} shows the error for the case where the number of training data is $200$. Moreover, error for the case where the number of training data is $120$ is demonstrated in Fig. \ref{n3}.

As we see, the results of performing the algorithm with the NICA initialization approach are better than the other methods. In fact, this strategy improved the convergence results. This method has the best results for accuracy and the number of iterations. As shown in the figure,  the NICA strategy has the best performance.  In fact, the algorithm based on this method converges faster, and the value of error is the lowest. In fact, as mentioned before, the NICA method is not suitable for some datasets due to the production of orthogonal matrices. These results also show that although the random strategy is not expensive,  its results are not as good as other methods.

	\section{Conclusion}
 In this paper, we studied various initialization strategies for NMF algorithms. In this direction,  we first reviewed some common approaches for solving NMF.  Then we calcified the initialization approaches for NMF. In fact, we divided these methods into four classes, that is, Random, Clastring, Heuristic, and Low-rank approximation methods. We summarize the advantages and disadvantages for some of them in Table \ref{conc}. As we see that, the random strategies, that is, random, random Acol, and Random C are cheap are methods. On the other hand, NNDSVD, FC-Means, Co-occurrence, and k-means are expensive methods. The NNSVD-LRC is another method for initialization that is not expensive than SVD. NICA and NPCA are two other methods that have are not expensive as the K-means strategy but they have better results. Numerical results presented in this paper showed that NICA has better results than the other strategies in terms of error. In fact, we performed  Lee's Algorithm on ORL datasets with three different cases, that is, the number of training data is $\{120,200,320\}$. The NICA has the lowest error and this strategy leads the algorithm to converge faster than the other strategies. 
	
	\begin{table}[h]
		\caption{Advantages and disadvantages for initialization methods}
		\centering
		{\color[rgb]{0.22,0.17,0.40}$\begin{array}{ *{10}{c} }
			\toprule
			\text{name} &\text{Pros} & \text{Cons} \\
			\midrule
			\cellcolor[gray]{0.8} \text{Random} & \text{easy}  & \text{dense $(W^0,H^0)$}    \\ 
			\cellcolor[gray]{0.8} & \text{cheap to compute}  & \text{with no intuitive meaning}   &  \\ \midrule
				\cellcolor[gray]{0.8}\text{Co-occurrence}  &   \text{uses term-term similarities}   &  \text{very expensive}  \\
			\bottomrule
			\cellcolor[gray]{0.8}\text{Random C }  &   \text{cheap }   &  \text{not very effective}        \\ 
			\bottomrule
			\cellcolor[gray]{0.8}\text{Random Aco}  &   \text{cheap }   &  \text{only slight decrease number of iterations}  \\
			\cellcolor[gray]{0.8} &   \text{sparse matrices built from original data}   &         \\ \midrule			
			\cellcolor[gray]{0.8} \text{k-means} &  \text{reduces NMF iterations}  & \text{dense}    \\ 
			\cellcolor[gray]{0.8}  &  \text{intuitive meaning of $W^0$}  &\text{expensive}  &\\ \midrule
			\cellcolor[gray]{0.8} \text{FC-means}&  \text{intuitive meaning of $(W^0,H^0)$}  & \text{dense $(W^0,H^0)$}  & \\ 
			\cellcolor[gray]{0.8} &      &   \text{expensive}     \\ \midrule
			\cellcolor[gray]{0.8}\text{NNDSVD}  &   \text{no randomization;}   &  \text{expensive}     \\
			\bottomrule
						\cellcolor[gray]{0.8}\text{NNSVD-LRC}  &   \text{no randomization;}   &  \text{expensive}     \\
				\cellcolor[gray]{0.8}  &   \text{non-negative;}   &    \\
			\bottomrule

			\bottomrule
			
		\end{array}$ }
	\label{conc}
	\end{table}

\subsection*{Abbreviations}
The following abbreviations are used in this paper:\\
NMF \qquad \qquad Nonnegative matrix factorizations\\
LR \qquad \qquad ~~ Low-rank\\
SVD  \qquad \qquad  Singular Value Decomposition\\
PCA  \qquad \qquad  Principal Component Analysis\\
KL  \qquad \qquad~~ Kullback-Leiber\\
MU \qquad \qquad ~  Mutilplicative Update\\
NPCA  \qquad \qquad Non-negative PCA\\
ICA \qquad \qquad  ~ Indipendent Component Analysis\\
NICA  \qquad \qquad  non-negative ICA\\
NNDSVD \qquad \qquad  Non-negative Double SVD\\
NNSVD-LRC \qquad  Non-negative SVD LR Correction\\
FCM\qquad \qquad Fuzzy C-means\\
DE \qquad \qquad  ~ Differential Evolution\\
PSO  \qquad \qquad Particle Swarm Optimization

\bibliographystyle{unsrtnat}
\bibliography{references}  

\begin{thebibliography}{83}
\providecommand{\natexlab}[1]{#1}
\providecommand{\url}[1]{\texttt{#1}}
\expandafter\ifx\csname urlstyle\endcsname\relax
  \providecommand{\doi}[1]{doi: #1}\else
  \providecommand{\doi}{doi: \begingroup \urlstyle{rm}\Url}\fi

\bibitem[Friedland et~al.(2011)Friedland, Mehrmann, Miedlar, and
  Nkengla]{Friedland}
Shmuel Friedland, V~Mehrmann, A~Miedlar, and M~Nkengla.
\newblock Fast low rank approximations of matrices and tensors.
\newblock \emph{The Electronic Journal of Linear Algebra}, 22:\penalty0
  1031--1048, 2011.

\bibitem[Eld{\'e}n(2003)]{Elden}
Lars Eld{\'e}n.
\newblock Numerical linear algebra and applications in data mining and it.
\newblock 2003.

\bibitem[Skillicorn(2007)]{Skillicorn}
David Skillicorn.
\newblock \emph{Understanding complex datasets: data mining with matrix
  decompositions}.
\newblock Chapman and Hall/CRC, 2007.

\bibitem[Murphy(2012)]{Murphy}
Kevin~P Murphy.
\newblock \emph{Machine learning: a probabilistic perspective}.
\newblock MIT press, 2012.

\bibitem[Lee et~al.(2013)Lee, Kim, Lebanon, and Singer]{Kim}
Joonseok Lee, Seungyeon Kim, Guy Lebanon, and Yoram Singer.
\newblock Matrix approximation under local low-rank assumption.
\newblock \emph{arXiv preprint arXiv:1301.3192}, 2013.

\bibitem[Golub and Van~Loan(1996)]{Golub}
Gene~H Golub and Charles~F Van~Loan.
\newblock Matrix computations 3rd edition.
\newblock \emph{The John Hopkins University, Baltimore}, 1996.

\bibitem[Alter et~al.(2000)Alter, Brown, and Botstein]{8}
Orly Alter, Patrick~O Brown, and David Botstein.
\newblock Singular value decomposition for genome-wide expression data
  processing and modeling.
\newblock \emph{Proceedings of the National Academy of Sciences}, 97\penalty0
  (18):\penalty0 10101--10106, 2000.

\bibitem[Wall et~al.(2003)Wall, Rechtsteiner, and Rocha]{17}
Michael~E Wall, Andreas Rechtsteiner, and Luis~M Rocha.
\newblock Singular value decomposition and principal component analysis.
\newblock In \emph{A practical approach to microarray data analysis}, pages
  91--109. Springer, 2003.

\bibitem[Datta(2010)]{Datta}
Biswa~Nath Datta.
\newblock \emph{Numerical linear algebra and applications}, volume 116.
\newblock Siam, 2010.

\bibitem[Sundarapandian(2008)]{Sundarapandian}
V~Sundarapandian.
\newblock \emph{Numerical linear algebra}.
\newblock PHI Learning Pvt. Ltd., 2008.

\bibitem[Trefethen and Bau~III(1997)]{Trefethen}
Lloyd~N Trefethen and David Bau~III.
\newblock \emph{Numerical linear algebra}, volume~50.
\newblock Siam, 1997.

\bibitem[Recht et~al.(2010)Recht, Fazel, and Parrilo]{Fazel}
Benjamin Recht, Maryam Fazel, and Pablo~A Parrilo.
\newblock Guaranteed minimum-rank solutions of linear matrix equations via
  nuclear norm minimization.
\newblock \emph{SIAM review}, 52\penalty0 (3):\penalty0 471--501, 2010.

\bibitem[Paruolo(2000)]{Reinsel}
Paolo Paruolo.
\newblock Multivariate reduced rank regression, theory and applications, 2000.

\bibitem[Miller and de~Callafon(2012)]{Miller}
Daniel~N Miller and Raymond~A de~Callafon.
\newblock Identification of linear time-invariant systems via constrained
  step-based realization.
\newblock \emph{IFAC Proceedings Volumes}, 45\penalty0 (16):\penalty0
  1155--1160, 2012.

\bibitem[Chu et~al.(2003)Chu, Funderlic, and Plemmons]{Chu}
Moody~T Chu, Robert~E Funderlic, and Robert~J Plemmons.
\newblock Structured low rank approximation.
\newblock \emph{Linear algebra and its applications}, 366:\penalty0 157--172,
  2003.

\bibitem[Drineas et~al.(2006)Drineas, Kannan, and Mahoney]{Petros}
Petros Drineas, Ravi Kannan, and Michael~W Mahoney.
\newblock Fast monte carlo algorithms for matrices iii: Computing a compressed
  approximate matrix decomposition.
\newblock \emph{SIAM Journal on Computing}, 36\penalty0 (1):\penalty0 184--206,
  2006.

\bibitem[Smilde et~al.(2005)Smilde, Bro, and Geladi]{Smilde}
Age Smilde, Rasmus Bro, and Paul Geladi.
\newblock \emph{Multi-way analysis: applications in the chemical sciences}.
\newblock John Wiley \& Sons, 2005.

\bibitem[Meyer(2000)]{Meyer}
Carl~D Meyer.
\newblock \emph{Matrix analysis and applied linear algebra}, volume~71.
\newblock Siam, 2000.

\bibitem[Kolda and O'leary(1998)]{Tamara}
Tamara~G Kolda and Dianne~P O'leary.
\newblock A semidiscrete matrix decomposition for latent semantic indexing
  information retrieval.
\newblock \emph{ACM Transactions on Information Systems (TOIS)}, 16\penalty0
  (4):\penalty0 322--346, 1998.

\bibitem[Esposito et~al.(2021)Esposito, Del~Buono, and Selicato]{23}
Flavia Esposito, Nicoletta Del~Buono, and Laura Selicato.
\newblock Nonnegative matrix factorization models for knowledge extraction from
  biomedical and other real world data.
\newblock \emph{PAMM}, 20\penalty0 (1):\penalty0 e202000032, 2021.

\bibitem[Meng et~al.(2016)Meng, Zeleznik, Thallinger, Kuster, Gholami, and
  Culhane]{14}
Chen Meng, Oana~A Zeleznik, Gerhard~G Thallinger, Bernhard Kuster, Amin~M
  Gholami, and Aed{\'\i}n~C Culhane.
\newblock Dimension reduction techniques for the integrative analysis of
  multi-omics data.
\newblock \emph{Briefings in bioinformatics}, 17\penalty0 (4):\penalty0
  628--641, 2016.

\bibitem[Sompairac et~al.(2019)Sompairac, Nazarov, Czerwinska, Cantini, Biton,
  Molkenov, Zhumadilov, Barillot, Radvanyi, Gorban, et~al.]{18}
Nicolas Sompairac, Petr~V Nazarov, Urszula Czerwinska, Laura Cantini, Anne
  Biton, Askhat Molkenov, Zhaxybay Zhumadilov, Emmanuel Barillot, Francois
  Radvanyi, Alexander Gorban, et~al.
\newblock Independent component analysis for unraveling the complexity of
  cancer omics datasets.
\newblock \emph{International Journal of molecular sciences}, 20\penalty0
  (18):\penalty0 4414, 2019.

\bibitem[Lee and Seung(1999)]{Daniel}
Daniel~D Lee and H~Sebastian Seung.
\newblock Learning the parts of objects by non-negative matrix factorization.
\newblock \emph{Nature}, 401\penalty0 (6755):\penalty0 788--791, 1999.

\bibitem[Aggarwal and Reddy(2014)]{Aggarwal}
Charu~C Aggarwal and Chandan~K Reddy.
\newblock Data clustering.
\newblock \emph{Algorithms and applications. Chapman\&Hall/CRC Data mining and
  Knowledge Discovery series, Londra}, 2014.

\bibitem[Lee and Seung()]{Lee01}
DD~Lee and HS~Seung.
\newblock Algorithms for non-negative matrix factorization. nips (2000).
\newblock \emph{Google Scholar}, pages 556--562.

\bibitem[Boutsidis and Gallopoulos(2008)]{Boutsidis}
Christos Boutsidis and Efstratios Gallopoulos.
\newblock Svd based initialization: A head start for nonnegative matrix
  factorization.
\newblock \emph{Pattern recognition}, 41\penalty0 (4):\penalty0 1350--1362,
  2008.

\bibitem[Casalino et~al.(2014)Casalino, Del~Buono, and Mencar]{Casalino}
Gabriella Casalino, Nicoletta Del~Buono, and Corrado Mencar.
\newblock Subtractive clustering for seeding non-negative matrix
  factorizations.
\newblock \emph{Information Sciences}, 257:\penalty0 369--387, 2014.

\bibitem[Esposito(2021)]{math}
Flavia Esposito.
\newblock A review on initialization methods for nonnegative matrix
  factorization: Towards omics data experiments.
\newblock \emph{Mathematics}, 9\penalty0 (9):\penalty0 1006, 2021.

\bibitem[Gong and Zhang(2012)]{g3}
Pinghua Gong and Changshui Zhang.
\newblock Efficient nonnegative matrix factorization via projected newton
  method.
\newblock \emph{Pattern Recognition}, 45\penalty0 (9):\penalty0 3557--3565,
  2012.

\bibitem[Bonettini et~al.(2008)Bonettini, Zanella, and Zanni]{g1}
Silvia Bonettini, Riccardo Zanella, and Luca Zanni.
\newblock A scaled gradient projection method for constrained image deblurring.
\newblock \emph{Inverse problems}, 25\penalty0 (1):\penalty0 015002, 2008.

\bibitem[Guan et~al.(2012)Guan, Tao, Luo, and Yuan]{g2}
Naiyang Guan, Dacheng Tao, Zhigang Luo, and Bo~Yuan.
\newblock Nenmf: An optimal gradient method for nonnegative matrix
  factorization.
\newblock \emph{IEEE Transactions on Signal Processing}, 60\penalty0
  (6):\penalty0 2882--2898, 2012.

\bibitem[Huang et~al.(2015)Huang, Liu, and Zhou]{g4}
Yakui Huang, Hongwei Liu, and Shuisheng Zhou.
\newblock Quadratic regularization projected barzilai--borwein method for
  nonnegative matrix factorization.
\newblock \emph{Data mining and knowledge discovery}, 29\penalty0 (6):\penalty0
  1665--1684, 2015.

\bibitem[Fathi-Hafshejani and Moaberfard(2020)]{g5}
S~Fathi-Hafshejani and Z~Moaberfard.
\newblock An interior-point algorithm for linearly constrained convex
  optimization based on kernel function and application in non-negative matrix
  factorization.
\newblock \emph{Optimization and Engineering}, 21\penalty0 (3):\penalty0
  1019--1051, 2020.

\bibitem[Albright et~al.(2006)Albright, Cox, Duling, Langville, and
  Meyer]{Albright}
Russell Albright, James Cox, David Duling, Amy~N Langville, and C~Meyer.
\newblock Algorithms, initializations, and convergence for the nonnegative
  matrix factorization.
\newblock Technical report, Tech. rep. 919. NCSU Technical Report Math 81706.
  http://meyer. math. ncsu~…, 2006.

\bibitem[Zheng et~al.(2007)Zheng, Yang, and Zhu]{Zheng07}
Zhonglong Zheng, Jie Yang, and Yitan Zhu.
\newblock Initialization enhancer for non-negative matrix factorization.
\newblock \emph{Engineering Applications of Artificial Intelligence},
  20\penalty0 (1):\penalty0 101--110, 2007.

\bibitem[Xue et~al.(2008)Xue, Tong, Chen, and Chen]{Xue08}
Yun Xue, Chong~Sze Tong, Ying Chen, and Wen-Sheng Chen.
\newblock Clustering-based initialization for non-negative matrix
  factorization.
\newblock \emph{Applied Mathematics and Computation}, 205\penalty0
  (2):\penalty0 525--536, 2008.

\bibitem[Alshabrawy et~al.(2012)Alshabrawy, Ghoneim, Awad, and
  Hassanien]{Alshabrawy}
Ossama~S Alshabrawy, ME~Ghoneim, WA~Awad, and Aboul~Ella Hassanien.
\newblock Underdetermined blind source separation based on fuzzy c-means and
  semi-nonnegative matrix factorization.
\newblock In \emph{2012 Federated Conference on Computer Science and
  Information Systems (FedCSIS)}, pages 695--700. IEEE, 2012.

\bibitem[Rezaei et~al.(2011)Rezaei, Boostani, and Rezaei]{Rezaei}
M~Rezaei, R~Boostani, and M~Rezaei.
\newblock An efficient initialization method for nonnegative matrix
  factorization.
\newblock \emph{Journal of Applied Sciences}, 11\penalty0 (2):\penalty0
  354--359, 2011.

\bibitem[Kim and Choi(2007)]{kim07}
Yong-Deok Kim and Seungjin Choi.
\newblock A method of initialization for nonnegative matrix factorization.
\newblock In \emph{2007 IEEE International Conference on Acoustics, Speech and
  Signal Processing-ICASSP'07}, volume~2, pages II--537. IEEE, 2007.

\bibitem[Stadlthanner et~al.(2007)Stadlthanner, Lutter, Theis, Lang, Tom{\'e},
  Georgieva, and Puntonet]{Stadlthanner}
Kurt Stadlthanner, Dominik Lutter, Fabian~J Theis, Elmar~Wolfgang Lang,
  Ana~Maria Tom{\'e}, Petia Georgieva, and Carlos~Garc{\'\i}a Puntonet.
\newblock Sparse nonnegative matrix factorization with genetic algorithms for
  microarray analysis.
\newblock In \emph{2007 International Joint Conference on Neural Networks},
  pages 294--299. IEEE, 2007.

\bibitem[Sn{\'a}{\v{s}}el et~al.(2008{\natexlab{a}})Sn{\'a}{\v{s}}el,
  Plato{\v{s}}, and Kr{\"o}mer]{Snasel}
V{\'a}clav Sn{\'a}{\v{s}}el, Jan Plato{\v{s}}, and Pavel Kr{\"o}mer.
\newblock Developing genetic algorithms for boolean matrix factorization.
\newblock \emph{Databases, Texts}, 61, 2008{\natexlab{a}}.

\bibitem[Price et~al.(2006)Price, Storn, and Lampinen]{Price}
Kenneth Price, Rainer~M Storn, and Jouni~A Lampinen.
\newblock \emph{Differential evolution: a practical approach to global
  optimization}.
\newblock Springer Science \& Business Media, 2006.

\bibitem[Atif et~al.(2019)Atif, Qazi, and Gillis]{syed}
Syed~Muhammad Atif, Sameer Qazi, and Nicolas Gillis.
\newblock Improved svd-based initialization for nonnegative matrix
  factorization using low-rank correction.
\newblock \emph{Pattern Recognition Letters}, 122:\penalty0 53--59, 2019.

\bibitem[Zhao et~al.(2008)Zhao, Zhuang, and Xu]{Zhao08}
Lihong Zhao, Guibin Zhuang, and Xinhe Xu.
\newblock Facial expression recognition based on pca and nmf.
\newblock In \emph{2008 7th World Congress on Intelligent Control and
  Automation}, pages 6826--6829. IEEE, 2008.

\bibitem[Kitamura and Ono(2016)]{Kitamura}
Daichi Kitamura and Nobutaka Ono.
\newblock Efficient initialization for nonnegative matrix factorization based
  on nonnegative independent component analysis.
\newblock In \emph{2016 IEEE International Workshop on Acoustic Signal
  Enhancement (IWAENC)}, pages 1--5. IEEE, 2016.

\bibitem[Oja and Plumbley(2004)]{Oja}
Erkki Oja and Mark Plumbley.
\newblock Blind separation of positive sources by globally convergent gradient
  search.
\newblock \emph{Neural Computation}, 16\penalty0 (9):\penalty0 1811--1825,
  2004.

\bibitem[Benachir et~al.(2013)Benachir, Hosseini, Deville, Karoui, and
  Hameurlain]{Benachir}
Djaouad Benachir, Shahram Hosseini, Yannick Deville, Moussa~Sofiane Karoui, and
  Abdelkader Hameurlain.
\newblock Modified independent component analysis for initializing non-negative
  matrix factorization: An approach to hyperspectral image unmixing.
\newblock In \emph{2013 IEEE 11th International Workshop of Electronics,
  Control, Measurement, Signals and their application to Mechatronics}, pages
  1--6. IEEE, 2013.

\bibitem[Sandler(2005)]{sandler}
Mark Sandler.
\newblock On the use of linear programming for unsupervised text
  classification.
\newblock In \emph{Proceedings of the eleventh ACM SIGKDD international
  conference on Knowledge discovery in data mining}, pages 256--264, 2005.

\bibitem[Liu and Tan(2018)]{pro1}
Zhaoqiang Liu and Vincent~YF Tan.
\newblock Rank-one nmf-based initialization for nmf and relative error bounds
  under a geometric assumption.
\newblock In \emph{2018 Information Theory and Applications Workshop (ITA)},
  pages 1--15. IEEE, 2018.

\bibitem[Zdunek(2012)]{pro2}
Rafal Zdunek.
\newblock Initialization of nonnegative matrix factorization with vertices of
  convex polytope.
\newblock In \emph{International Conference on Artificial Intelligence and Soft
  Computing}, pages 448--455. Springer, 2012.

\bibitem[Sauwen et~al.(2017)Sauwen, Acou, Bharath, Sima, Veraart, Maes,
  Himmelreich, Achten, and Van~Huffel]{pro3}
Nicolas Sauwen, Marjan Acou, Halandur~N Bharath, Diana~M Sima, Jelle Veraart,
  Frederik Maes, Uwe Himmelreich, Eric Achten, and Sabine Van~Huffel.
\newblock The successive projection algorithm as an initialization method for
  brain tumor segmentation using non-negative matrix factorization.
\newblock \emph{Plos one}, 12\penalty0 (8):\penalty0 e0180268, 2017.

\bibitem[Mahoney(2011)]{Mahoney}
Michael~W Mahoney.
\newblock Randomized algorithms for matrices and data.
\newblock \emph{arXiv preprint arXiv:1104.5557}, 2011.

\bibitem[Drineas and Mahoney(2016)]{Drineas}
Petros Drineas and Michael~W Mahoney.
\newblock Randnla: randomized numerical linear algebra.
\newblock \emph{Communications of the ACM}, 59\penalty0 (6):\penalty0 80--90,
  2016.

\bibitem[Wang and Li(2010)]{Wangli}
Fei Wang and Ping Li.
\newblock Efficient nonnegative matrix factorization with random projections.
\newblock In \emph{Proceedings of the 2010 SIAM International Conference on
  Data Mining}, pages 281--292. SIAM, 2010.

\bibitem[Tepper and Sapiro(2016)]{Tepper}
Mariano Tepper and Guillermo Sapiro.
\newblock Compressed nonnegative matrix factorization is fast and accurate.
\newblock \emph{IEEE Transactions on Signal Processing}, 64\penalty0
  (9):\penalty0 2269--2283, 2016.

\bibitem[Langville et~al.(2006)Langville, Meyer, Albright, Cox, and
  Duling]{Amy2}
Amy~N Langville, Carl~D Meyer, Russell Albright, James Cox, and David Duling.
\newblock Alternating least squares algorithms for the nonnegative matrix
  factorization, 2006.

\bibitem[Berthold et~al.(2010)Berthold, Borgelt, H{\"o}ppner, and
  Klawonn]{Berthold}
Michael~R Berthold, Christian Borgelt, Frank H{\"o}ppner, and Frank Klawonn.
\newblock \emph{Guide to intelligent data analysis: how to intelligently make
  sense of real data}.
\newblock Springer Science \& Business Media, 2010.

\bibitem[Xu et~al.(2003)Xu, Liu, and Gong]{Xu}
Wei Xu, Xin Liu, and Yihong Gong.
\newblock Document clustering based on non-negative matrix factorization.
\newblock In \emph{Proceedings of the 26th annual international ACM SIGIR
  conference on Research and development in informaion retrieval}, pages
  267--273, 2003.

\bibitem[Chr{\'e}tien et~al.(2016)Chr{\'e}tien, Guyeux, Conesa, Delage-Mouroux,
  Jouvenot, Huetz, and Desc{\^o}tes]{Chretien}
St{\'e}phane Chr{\'e}tien, Christophe Guyeux, Bastien Conesa, R{\'e}gis
  Delage-Mouroux, Mich{\`e}le Jouvenot, Philippe Huetz, and Fran{\c{c}}oise
  Desc{\^o}tes.
\newblock A bregman-proximal point algorithm for robust non-negative matrix
  factorization with possible missing values and outliers-application to gene
  expression analysis.
\newblock \emph{BMC bioinformatics}, 17\penalty0 (8):\penalty0 623--631, 2016.

\bibitem[Cichocki and Amari(2002)]{Cichocki4}
Andrzej Cichocki and Shun-ichi Amari.
\newblock \emph{Adaptive blind signal and image processing: learning algorithms
  and applications}.
\newblock John Wiley \& Sons, 2002.

\bibitem[Perronnin and Bouchard(2017)]{Shahnaz}
Florent Perronnin and Guillaume Bouchard.
\newblock Clustering using non-negative matrix factorization on sparse graphs,
  August~8 2017.
\newblock US Patent 9,727,532.

\bibitem[Wild et~al.(2004)Wild, Curry, and Dougherty]{Wild04}
Stefan Wild, James Curry, and Anne Dougherty.
\newblock Improving non-negative matrix factorizations through structured
  initialization.
\newblock \emph{Pattern recognition}, 37\penalty0 (11):\penalty0 2217--2232,
  2004.

\bibitem[Wild et~al.(2003)Wild, Wild, Curry, Dougherty, and Betterton]{Wild02}
Stefan Wild, Written~Stefan Wild, James Curry, Anne Dougherty, and Meredith
  Betterton.
\newblock \emph{Seeding non-negative matrix factorizations with the spherical
  k-means clustering}.
\newblock PhD thesis, Citeseer, 2003.

\bibitem[Forgey(1965)]{Forgy}
Edward Forgey.
\newblock Cluster analysis of multivariate data: Efficiency vs.
  interpretability of classification.
\newblock \emph{Biometrics}, 21\penalty0 (3):\penalty0 768--769, 1965.

\bibitem[Hartigan(1975)]{Hartigan}
John~A Hartigan.
\newblock \emph{Clustering algorithms}.
\newblock John Wiley \& Sons, Inc., 1975.

\bibitem[Din()]{Ding}


\bibitem[Anderberg(2014)]{Anderberg}
Michael~R Anderberg.
\newblock \emph{Cluster analysis for applications: probability and mathematical
  statistics: a series of monographs and textbooks}, volume~19.
\newblock Academic press, 2014.

\bibitem[Dong et~al.(2014)Dong, Lin, and Chu]{Dong}
Bo~Dong, Matthew~M Lin, and Moody~T Chu.
\newblock Nonnegative rank factorization—a heuristic approach via rank
  reduction.
\newblock \emph{Numerical Algorithms}, 65\penalty0 (2):\penalty0 251--274,
  2014.

\bibitem[Goldberg and Holland(1988)]{Goldberg}
David~E Goldberg and John~Henry Holland.
\newblock Genetic algorithms and machine learning.
\newblock 1988.

\bibitem[Sn{\'a}{\v{s}}el et~al.(2008{\natexlab{b}})Sn{\'a}{\v{s}}el,
  Plato{\v{s}}, and Kr{\"o}mer]{asel}
V{\'a}clav Sn{\'a}{\v{s}}el, Jan Plato{\v{s}}, and Pavel Kr{\"o}mer.
\newblock Developing genetic algorithms for boolean matrix factorization.
\newblock \emph{Databases, Texts}, 61, 2008{\natexlab{b}}.

\bibitem[Eberhart and Kennedy(1995)]{Kennedy}
Russell Eberhart and James Kennedy.
\newblock Particle swarm optimization.
\newblock In \emph{Proceedings of the IEEE international conference on neural
  networks}, volume~4, pages 1942--1948. Citeseer, 1995.

\bibitem[Janecek and Tan(2011)]{Janecek}
Andreas Janecek and Ying Tan.
\newblock Using population based algorithms for initializing nonnegative matrix
  factorization.
\newblock In \emph{International Conference in Swarm Intelligence}, pages
  307--316. Springer, 2011.

\bibitem[Deerwester et~al.(1990)Deerwester, Dumais, Furnas, Landauer, and
  Harshman]{S.Deerwester}
Scott Deerwester, Susan~T Dumais, George~W Furnas, Thomas~K Landauer, and
  Richard Harshman.
\newblock Indexing by latent semantic analysis.
\newblock \emph{Journal of the American society for information science},
  41\penalty0 (6):\penalty0 391--407, 1990.

\bibitem[Eckart and Young(1936)]{Eckart}
Carl Eckart and Gale Young.
\newblock The approximation of one matrix by another of lower rank.
\newblock \emph{Psychometrika}, 1\penalty0 (3):\penalty0 211--218, 1936.

\bibitem[Cao(2006)]{Lijie}
Lijie Cao.
\newblock Singular value decomposition applied to digital image processing.
\newblock \emph{Division of Computing Studies, Arizona State University
  Polytechnic Campus, Mesa, Arizona State University polytechnic Campus}, pages
  1--15, 2006.

\bibitem[Cichocki et~al.(2009)Cichocki, Zdunek, Phan, and Amari]{Cichocki3}
Andrzej Cichocki, Rafal Zdunek, Anh~Huy Phan, and Shun-ichi Amari.
\newblock \emph{Nonnegative matrix and tensor factorizations: applications to
  exploratory multi-way data analysis and blind source separation}.
\newblock John Wiley \& Sons, 2009.

\bibitem[Guillamet and Vitria(2003)]{Guillamet}
David Guillamet and Jordi Vitria.
\newblock Evaluation of distance metrics for recognition based on non-negative
  matrix factorization.
\newblock \emph{Pattern Recognition Letters}, 24\penalty0 (9-10):\penalty0
  1599--1605, 2003.

\bibitem[Turk and Pentland(1991)]{Turk}
Matthew Turk and Alex Pentland.
\newblock Eigenfaces for recognition.
\newblock \emph{Journal of cognitive neuroscience}, 3\penalty0 (1):\penalty0
  71--86, 1991.

\bibitem[Geng et~al.(2016)Geng, Ji, and Sun]{Xiu-rui-16}
Xiu-rui Geng, Lu-yan Ji, and Kang Sun.
\newblock Non-negative matrix factorization based unmixing for principal
  component transformed hyperspectral data.
\newblock \emph{Frontiers of Information Technology \& Electronic Engineering},
  17\penalty0 (5):\penalty0 403--412, 2016.

\bibitem[Shen and Huang(2008)]{Shen}
Haipeng Shen and Jianhua~Z Huang.
\newblock Sparse principal component analysis via regularized low rank matrix
  approximation.
\newblock \emph{Journal of multivariate analysis}, 99\penalty0 (6):\penalty0
  1015--1034, 2008.

\bibitem[Zass and Shashua(2007)]{Zass}
Ron Zass and Amnon Shashua.
\newblock Nonnegative sparse pca.
\newblock In \emph{Advances in neural information processing systems}, pages
  1561--1568. Citeseer, 2007.

\bibitem[Bauckhage(2014)]{29}
Christian Bauckhage.
\newblock A purely geometric approach to non-negative matrix factorization.
\newblock In \emph{LWA}, pages 125--136. Citeseer, 2014.

\bibitem[Smaragdis et~al.(2006)Smaragdis, Raj, and Shashanka]{30}
Paris Smaragdis, Bhiksha Raj, and Madhusudana Shashanka.
\newblock A probabilistic latent variable model for acoustic modeling.
\newblock \emph{Advances in models for acoustic processing, NIPS},
  148:\penalty0 8--1, 2006.

\end{thebibliography}






\end{document}